\documentclass[12pt]{article}
\usepackage{amssymb}
\usepackage{graphicx}
\usepackage{amsmath}

\newtheorem{theorem}{Theorem}

\newtheorem{lemma}[theorem]{Lemma}

\newtheorem{question}[theorem]{Question}

\begin{document}

\title{Less than continuum many translates of a compact nullset may cover any
infinite profinite group}
\author{Mikl\'{o}s Ab\'{e}rt}
\maketitle

\begin{abstract}
We show that it is consistent with the axioms of set theory that every
infinite profinite group $G$ possesses a closed subset $X$ of Haar measure
zero such that less than continuum many translates of $X$ cover $G$. This
answers a question of Elekes and T\'{o}th and by their work settles the
problem for all infinite compact topological groups.
\end{abstract}

Let $G$ be a profinite group and let $\mu $ be a left Haar measure on $G$.
If $X\subseteq G$ is a measurable subset of measure zero, then trivially
countably many $G$-translates of $X$ can not cover $G$. Thus the
continuum-hypothesis (CH) implies that one needs at least continuum many
translates of $X$ to do this. However, it turns out that if we do not assume
CH, then, consistently with the axioms of set theory, the world may be
different.

\begin{theorem}
\label{eredmeny}It is consistent with the axioms of set theory that for
every infinite profinite group $G$ there exists a closed subset $X\subseteq
G $ of Haar-measure $0$ such that less than continuum many $G$-translates of 
$X $ cover $G$.
\end{theorem}

Note that our proof is not constructive in the sense that it uses a
probabilistic argument.

Theorem \ref{eredmeny} answers a question of Elekes and T\'{o}th 
\cite[Question 3.10]{ET}. The first result of this type is due to Elekes and
Steprans \cite{ES} who proved the analogue of Theorem \ref{eredmeny} for the
real line. For the background of covering the real line with less than
continuum many nullsets see \cite{Bart}. Recently Elekes and T\'{o}th \cite
{ET} proved the theorem in the realm of real Lie groups, locally compact
Abelian groups and also showed that Theorem \ref{eredmeny} would imply the
result for arbitrary infinite compact topological groups. So now this is
settled.

In turn, Darji and T. Keleti \cite{DK} proved that if $X$ is a compact
subset of the reals with packing dimension less than $1$ then one needs
continuum many translates of $X$ to cover the real line. Let $(N_{i})$ be a
strictly descending chain of open normal subgroups in $G$. Then $(N_{i})$
naturally gives rise to a dimension function; for a subset $X\subseteq G$
let 
\begin{equation*}
\dim _{(N_{i})}(X)=\lim \inf_{i}\frac{\log \left| XN_{i}/N_{i}\right| }{\log
\left| G/N_{i}\right| }
\end{equation*}
See \cite{BS} for an introduction on this dimension notion. This suggests
the following.

\begin{question}
Let $G$ be an infinite profinite group and let $X\subseteq G$ be a compact
subset such that $\dim X<1$ with respect to some strictly descending chain
of open normal subgroups of $G$. Is it true that one needs at least
continuum many translates of $X$ to cover $G$?
\end{question}

Let $G$ be a group. We say that a subset $Y\subseteq G$ \emph{can be
translated into} $X\subseteq G$ if there exists $g\in G$ such that $%
gY\subseteq X$.

Let $G$ be a profinite group and let $N_{0},N_{1},\ldots $ be open normal
subgroups of $G$. Let $f:\mathbb{N}\rightarrow \mathbb{N}$ be a function. We
say that a subset $X\subseteq G$ is $f$\emph{-thin} with respect to $(N_{i})$%
, if 
\begin{equation*}
\left| XN_{i}/N_{i}\right| \leq f(i)\text{ (}i\in \mathbb{N}\text{).}
\end{equation*}
The proof of Theorem \ref{eredmeny} relies on the following result which is
purely group theoretical.

\begin{theorem}
\label{main}Let $G$ be an infinite profinite group. Let $f(0)=1$ and $f(n)=n$
($n>0$). Then there exists a descending chain of open normal subgroups in $%
G=N_{0}>N_{1}>N_{2}>\ldots $ and a closed subset $X\subseteq G$ of measure
zero such that every subset of $G$ that is $f$-thin with respect to $(N_{i})$
can be translated into $X$.
\end{theorem}

We prove this theorem after a series of results on finite groups.

\begin{theorem}
\label{gazdag}Let $G$ be a finite group of size $n$ and let 
\begin{equation*}
k<\frac{n-\log 2}{\log n}
\end{equation*}
be a natural number. Then for all integers $l$ with 
\begin{equation*}
2\left( k\log n+\log 2\right) ^{\frac{1}{k}}n^{1-\frac{1}{k}}<l\leq n
\end{equation*}
there exist subsets $X_{1},\ldots ,X_{k}\subseteq G$ of size $l$ such that
for every $g_{1},\ldots ,g_{k}\in G$ the intersection $\cap _{i}X_{i}g_{i}$
is non-empty.
\end{theorem}

\noindent \textbf{Proof. }We can trivially assume $n\geq 3$ and $k\geq 2$.

Let 
\begin{equation*}
p=\left( \frac{k\log n+\log 2}{n}\right) ^{\frac{1}{k}}\text{.}
\end{equation*}
The assumption on $k$ implies $p<1$. By a random subset of $G$ we mean a
subset which is produced by choosing each element of $G$ independently with
probability $p$. Let $X_{1},X_{2},\ldots ,X_{k}$ be independent random
subsets of $G$. We will show that with positive probability these subsets
satisfy the theorem.

Fix $g_{1},g_{2},\ldots ,g_{k}\in G$. Then the subsets $%
X_{1}g_{1},X_{2}g_{2},\ldots ,X_{k}g_{k}$ are also independent random
subsets of $G$. Thus for all $x\in X$ the probability 
\begin{equation*}
P(x\in \bigcap_{i=1}^{k}X_{i}g_{i})=\prod_{i=1}^{k}P(x\in X_{i}g_{i})=p^{k}
\end{equation*}
Also, the events 
\begin{equation*}
(x\in \bigcap_{i=1}^{k}X_{i}g_{i})\text{ \ }(x\in G)
\end{equation*}
are independent, so we have 
\begin{equation*}
P(\bigcap_{i=1}^{k}X_{i}g_{i}=\emptyset )=P(\text{for all }x\in G\text{ we
have }x\notin \bigcap_{i=1}^{k}X_{i}g_{i})=(1-p^{k})^{n}
\end{equation*}

Summing up these bad probabilities for the possible $g_{1},g_{2},\ldots
,g_{k}\in G$ we get 
\begin{equation*}
P=P(\text{there exists }g_{1},\ldots ,g_{k}\in G\text{ with }%
\bigcap_{i=1}^{k}X_{i}g_{i}=\emptyset )\leq n^{k}(1-p^{k})^{n}
\end{equation*}
and using 
\begin{equation*}
-\log (1-p^{k})>p^{k}=\frac{k\log n+\log 2}{n}
\end{equation*}
we get 
\begin{equation*}
P\leq n^{k}(1-p^{k})^{n}<\frac{1}{2}
\end{equation*}

We also want to control the sizes of the $X_{i}$ from above. For each $1\leq
i\leq k$ the size $\left| X_{i}\right| $ equals the sum of independent
variables which take $1$ with probability $p$ and $0$ with probability $1-p$%
. Thus the expected value is 
\begin{equation*}
E(\left| X_{i}\right| )=pn=\frac{1}{2}n^{1-\frac{1}{k}}(k-1)\log n
\end{equation*}
and the standard deviance is 
\begin{equation*}
s=SD(\left| X_{i}\right| )=\sqrt{p(1-p)n}<\sqrt{pn}
\end{equation*}
Now using Chernoff's inequality 
\begin{equation*}
P(\left| X_{i}\right| >as)\leq 2\exp (-a^{2}/4)
\end{equation*}
for $a=2pn/s$ we obtain 
\begin{equation*}
Q_{i}=P(\left| X_{i}\right| >2pn)\leq 2\exp (-\frac{p^{2}n^{2}}{s^{2}}%
)<2\exp (-pn)
\end{equation*}

Using $n\geq 3$ and $k\geq 2$ we get 
\begin{equation*}
k<\frac{n-\log 2}{\log n}<\frac{n^{2}}{4}
\end{equation*}
which implies 
\begin{equation*}
p=\left( \frac{k\log n+\log 2}{n}\right) ^{\frac{1}{k}}>\frac{2\log n}{n}>%
\frac{\log 4k}{n}
\end{equation*}
yielding 
\begin{equation*}
Q_{i}<2\exp (-pn)<2\exp (-\log 4k)=\frac{1}{2k}
\end{equation*}
Adding up, we obtain 
\begin{equation*}
P+\sum_{i=1}^{k}Q_{i}<1
\end{equation*}
which implies that with positive probability for all $1\leq i\leq k$ we have 
\begin{equation*}
\left| X_{i}\right| \leq 2pn=2\left( k\log n+\log 2\right) ^{\frac{1}{k}%
}n^{1-\frac{1}{k}}
\end{equation*}
and for all $g_{1},g_{2},\ldots ,g_{k}\in G$ we have 
\begin{equation*}
\bigcap_{i=1}^{k}X_{i}g_{i}\neq \emptyset \text{.}
\end{equation*}
In particular, these events happen for at least one value of $%
X_{1},X_{2},\ldots ,X_{k}$. If necessary, we can enlarge the $X_{i}$ to have
size $l$, keeping the intersection condition untouched. The theorem is
proved. $\square $

\bigskip

Let $k$ be a natural number. We say that a subset $X\subseteq G$ is $k$\emph{%
-covering} in $G$ if every subset $Y\subseteq G$ of size $k$ can be
translated into $X$.

\begin{lemma}
\label{perketto}Let $G$ be a finite group of size $n$ and let $k$ be a
natural number. If 
\begin{equation*}
(4k)^{k}\left( k\log n+\log 2\right) <n
\end{equation*}
then there exists a $k$-covering subset $X\subseteq G$ of size $\left|
X\right| \leq n/2$.
\end{lemma}

\noindent \textbf{Proof. }The inequality implies 
\begin{equation*}
2\left( k\log n+\log 2\right) ^{\frac{1}{k}}n^{1-\frac{1}{k}}<\frac{n}{2k}
\end{equation*}
Since the assumptions of Theorem \ref{gazdag} hold, there exist subsets $%
X_{1},X_{2},\ldots ,X_{k}\subseteq G$ of size at most $n/2k$ satisfying the
conclusion. Let $X=\cup _{i}X_{i}$. Then we have $\left| X\right| \leq n/2$.

Let $Y\subseteq G$ be an arbitrary subset of size $k$. Let us list the
elements of $Y$ as $y_{1},y_{2},\ldots ,y_{k}$ and let 
\begin{equation*}
K=\cap _{i}X_{i}y_{i}^{-1}\text{.}
\end{equation*}
By Theorem \ref{gazdag} $K$ is non-empty. Let $g\in K$. Now for $1\leq i\leq
k$ we have 
\begin{equation*}
gy_{i}\in X_{i}y_{i}^{-1}y_{i}=X_{i}
\end{equation*}
implying 
\begin{equation*}
gY\subseteq \cup _{i}X_{i}=X.
\end{equation*}
So $X$ is a $k$-covering subset of $G$. $\square $

\bigskip

The following extension lemma is the key to pass from finite to profinite
groups.

\begin{lemma}
\label{ragaszt}Let $G$ be a finite group, let $\varphi :G\rightarrow H$ be
an epimorphism with kernel $N$ of size $n$ and let $X\subseteq H$ be a
subset in $H$. Assume that 
\begin{equation*}
(4k)^{k}\left( k\log n+\log 2\right) <n
\end{equation*}
Then there exists a subset $X^{\prime }\subseteq G$ such that \newline
1) $\varphi (X^{\prime })=X$;\newline
2) $\left| X^{\prime }\right| \leq n\left| X\right| /2$;\newline
3) Every subset $Y\subseteq G$ of size at most $k+1$ for which $\varphi (Y)$
can be translated into $X$ can be translated into $X^{\prime }$.
\end{lemma}

\noindent \textbf{Proof. }Let $\widetilde{X}\subseteq G$ such that $\varphi (%
\widetilde{X})=X$ and $\left| \widetilde{X}\right| =\left| X\right| $. Lemma 
\ref{perketto} implies that there exists a $k+1$-covering subset $L\subseteq
N$ of size $\left| L\right| \leq n/2$. Let $X^{\prime }=L\widetilde{X}$.
Then 1) and 2) trivially hold and we only have to check 3).

Let $Y$ $\subseteq G$ be a subset with for which $\varphi (Y)$ can be
translated into $X$. This means that there exists $h\in H$ such that $%
h\varphi (Y)\subseteq X$. Let $g\in G$ such that $\varphi (g)=h$. Then we
have $\varphi (gY)=h\varphi (Y)\subseteq X$, implying $gY\subseteq N%
\widetilde{X}$. Let us list the elements of $\widetilde{X}$ as $%
x_{1},x_{2},\ldots ,x_{t}$ with $t=\left| X\right| $. Using this notation,
there exist subsets $Y_{1},Y_{2},\ldots ,Y_{t}\subseteq N$ such that 
\begin{equation*}
gY=\bigcup_{i=1}^{t}Y_{i}x_{i}\text{.}
\end{equation*}

Let $Y^{\prime }=\cup _{i}Y_{i}\subseteq N$. Then $\left| Y^{\prime }\right|
\leq \left| Y\right| =k+1$, so there exists an element $u\in N$ such that $%
uY^{\prime }\subseteq L$. Now 
\begin{equation*}
ugY=\bigcup_{i=1}^{t}uY_{i}x_{i}\subseteq \bigcup_{i=1}^{t}Lx_{i}=L%
\widetilde{X}=X^{\prime }
\end{equation*}
and so $Y$ can be translated into $X^{\prime }$. $\square $

\bigskip

\noindent \textbf{Proof of Theorem \ref{main}. }Since $G$ is infinite, there
exists a descending chain of open normal subgroups $G=N_{0}>N_{1}>N_{2}>%
\ldots $ Let $n_{i}=\left| N_{i}/N_{i+1}\right| $. By passing to a suitable
subchain we can also assume that for all $i>0$ we have 
\begin{equation*}
(4i)^{i}\left( i\log n_{i-1}+\log 2\right) <n_{i-1}
\end{equation*}
Let $G_{i}=G/N_{i}$ and let $\varphi _{i}:G\rightarrow G_{i}$ denote the
quotient map.

We claim that there exist subsets $X_{i}\in G_{i}$ such that for all $i\in 
\mathbb{N}$ the following hold:

\begin{enumerate}
\item  $\varphi _{i}(X_{i+1}N_{i+1})\subseteq X_{i}$;

\item  $\left| X_{i}\right| \leq \left| G_{i}\right| /2^{i}$;

\item  If $Y\subseteq G$ is an $f$-thin subset with respect to $(N_{i})$
then $\varphi _{i}(Y)$ can be translated into $X_{i}$ in $G_{i}$.
\end{enumerate}

We prove our claims by induction on $i$. For $i=0$ we set $X_{i}=G_{i}=\{1\}$
and everything holds trivially.

Assume that we already constructed $X_{i}$ as above. Let $\varphi
:G_{i+1}\rightarrow G_{i}$ be defined by $\varphi (g)=\varphi _{i}(gN_{i})$
and let $H=G_{i}$. Then setting $X=X_{i}$ and $G=G_{i+1}$ we can use Lemma 
\ref{ragaszt}. Let $X_{i+1}=X^{\prime }$ obtained from the lemma. Now 1)
follows immediately from 1) of Lemma \ref{ragaszt}. By 2) of Lemma \ref
{ragaszt} and using induction we have 
\begin{equation*}
\left| X_{i+1}\right| \leq \left| X_{i}\right| \left| \ker \varphi \right|
/2\leq \left| G_{i}\right| \left| \ker \varphi \right| /2^{i+1}=\left|
G_{i+1}\right| /2^{i+1}
\end{equation*}
so 2) also holds. Finally, if $Y\subseteq G$ is an $f$-thin subset with
respect to $(N_{i})$ then $\varphi _{i}(Y)\subseteq G_{i}$ can be translated
into $X_{i}$ by induction. Since $Y\subseteq G$ is $f$-thin, we have 
\begin{equation*}
\left| \varphi _{i+1}(YN_{i+1})\right| \leq i+1
\end{equation*}
and so 3) of Lemma \ref{ragaszt} holds for $\varphi _{i+1}(Y)$. It follows
that $\varphi _{i+1}(Y)$ can be translated into $X_{i+1}$ in $G_{i+1}$ what
we claimed in 3). So all the claims hold.

Now let 
\begin{equation*}
X=\bigcap_{i}X_{i}N_{i}\subseteq G\text{.}
\end{equation*}
Since all the $X_{i}N_{i}$ are open and closed, $X$ is closed. Also, if $\mu 
$ stands for the normalised Haar measure on $G$ then using 2) we have 
\begin{equation*}
\mu (X)\leq \mu (X_{i}N_{i})=\frac{\left| X_{i}\right| }{\left| G_{i}\right| 
}\leq \frac{1}{2^{i}}\text{ (}i\in \mathbb{N}\text{)}
\end{equation*}
which implies $\mu (X)=0$.

Let $Y\subseteq G$ be an $f$-thin subset of $G$ with respect to $(N_{i})$.
For $i\in \mathbb{N}$ let 
\begin{equation*}
T_{i}=\left\{ g\in G\mid \varphi _{i}(gY)\subseteq X_{i}\right\} \text{.}
\end{equation*}
Then 3) tells us that the $T_{i}$ are non-empty. Also the $T_{i}$ are unions
of $N_{i}$-cosets, so they are open and closed. On the other hand, if $g\in
T_{i+1}$ then $\varphi _{i+1}(gY)\subseteq X_{i+1}$ and so, using 1) we have 
\begin{equation*}
\varphi _{i}(gY)=\varphi _{i}(\varphi _{i+1}(gY)N_{i+1})\subseteq \varphi
_{i}(X_{i+1}N_{i+1})\subseteq X_{i}
\end{equation*}
which implies $g\in T_{i}$. That is, we have $T_{i+1}\subseteq T_{i}$ ($i\in 
\mathbb{N}$). Now by the compactness of $G$ the $T_{i}$ have non-empty
intersection. Let $g\in \cap _{i}T_{i}$.

From $\varphi _{i}(gY)\subseteq X_{i}$ we have $gY\subseteq X_{i}N_{i}$.
Thus $gY\subseteq \cap _{i}X_{i}N_{i}=X$ and so $Y$ can be translated into $%
X $. $\square $

\bigskip

We are ready to prove Theorem \ref{eredmeny}.

\bigskip

\noindent \textbf{Proof of Theorem \ref{eredmeny}. }Let $f(0)=1$ and $f(n)=n$
($n>0$). Applying Theorem \ref{main} to $G$ we get a descending chain of
open normal subgroups in $G=N_{0}>N_{1}>N_{2}>\ldots $ and a closed subset $%
X\subseteq G$ of measure zero such that every subset of $G$ that is $f$-thin
with respect to $(N_{i})$ can be translated into $X$. For $i\in \mathbb{N}$
let $G_{i}=G/N_{i}$ and let $\varphi _{i}:G\rightarrow G_{i}$ denote the
quotient map.

Let $P=\prod_{i}G_{i}$ and let $\varphi :G\rightarrow P$ be defined by 
\begin{equation*}
\varphi (g)=(\varphi _{0}(g),\varphi _{1}(g),\ldots ,\varphi _{n}(g),\ldots )
\end{equation*}
Then $P$ is a profinite group endowed with the product topology, $\varphi $
is continuous and since $G$ is compact, $A=\varphi (G)$ is a closed subgroup
of $P$.

Following the notation in \cite{ET} we call a subset $S=\prod_{i}S_{i}%
\subseteq \prod_{i}G_{i}=P$ an $\dot{f}$\emph{-slalom} if $\left|
S_{i}\right| \leq f(i)$ ($i\in \mathbb{N}$). It is consistent with the
axioms of set theory (see \cite{Bart}) that there exists a cardinality $%
\kappa <2^{\varpi }$ such that we can cover $P$ by $\kappa $ $f$-slaloms.
Let $L_{j}$ ($j\in J$) be a set of $f$-slaloms such that $\cup _{j\in
J}L_{j}=P$. Here $J$ is an index set of cardinality $\kappa $. For $j\in J$
let $M_{j}=A\cap L_{j}$ and let $Y_{j}=\varphi ^{-1}(M_{j})$. Then $A=\cup
_{j\in J}M_{j}$ and so $G=\cup _{j\in J}Y_{j}$.

Now for all $j\in Y$ the subset $Y_{j}$ is $f$-thin with respect to $(N_{i})$
since $\varphi (Y_{j})=M_{j}\subseteq L_{j}$. So by Theorem \ref{main} there
exists $g_{j}\in G$ such that $g_{j}Y_{j}\subseteq X$, that is, $%
Y_{j}\subseteq g_{j}^{-1}X$. But then 
\begin{equation*}
\bigcup_{j\in J}g_{j}^{-1}X\supseteq \bigcup_{j\in J}Y_{j}=G
\end{equation*}
and so we covered $G$ by $\kappa <2^{\varpi }$ translates of $X$. $\square $

\bigskip

\noindent \textbf{Remark. }For a finite group $G$ of size $n$ and a natural
number $k$ let $\mathrm{cov}(G,k)$ denote the size of a minimal $k$-covering
subset of $G$. Using Theorem \ref{gazdag} and the proof of Lemma \ref
{perketto} we get 
\begin{equation*}
\mathrm{cov}(G,k)\leq 2k\left( k\log n+\log 2\right) ^{\frac{1}{k}}n^{1-%
\frac{1}{k}}
\end{equation*}
On the other hand, it is easy to show that if $X$ is a $k$-covering subset
of $G$ then for all $g_{1},g_{2},\ldots ,g_{k}\in G$ we have $\cap
_{i}X_{i}g_{i}\neq \emptyset $. We claim that this implies $\left| X\right|
\geq n^{1-\frac{1}{k}}$. Indeed, if $A$ and $B$ are subsets of $G$, then it
is easy to see that the expected value of $\left| A\cap Bg\right| =\left|
A\right| \left| B\right| /n$ where $g$ is a uniform random element of $G$.
So with the notation $g_{1}=1$ and $l=\left| X\right| $ there exists $%
g_{2}\in G$ such that $\left| Xg_{1}\cap Xg_{2}\right| <l^{2}/n$. Similarly,
we can find $g_{3}\in G$ such that $\left| Xg_{1}\cap Xg_{2}\cap
Xg_{3}\right| <l^{3}/n^{2}$. Using induction, we can find $g_{4},\ldots
,g_{k}\in G$ such that $\left| \cap _{i}X_{i}g_{i}\right| <l^{n}/k^{n-1}$.
Now if $\left| X\right| <n^{1-\frac{1}{k}}$ then $\left| \cap
_{i}X_{i}g_{i}\right| <1$ that is, it is empty. This shows that our claim
holds and that 
\begin{equation*}
n^{1-\frac{1}{k}}\leq \mathrm{cov}(G,k)
\end{equation*}
It would be interesting to see whether the extra $2k\left( k\log n+\log
2\right) ^{\frac{1}{k}}$ in the upper estimate is really needed. For the
case $k=2$ the question is equivalent to finding the smallest subset $X$ of $%
G$ with $XX^{-1}=G$. This has been settled in \cite{fkl} (see also \cite{KL}%
). It follows that 
\begin{equation*}
\mathrm{cov}(G,2)=O(\sqrt{n})
\end{equation*}
so this suggests that at least the $\log n$ factor in $\mathrm{cov}(G,k)$
could be omitted. Note, however, that elementary counting arguments only
give $\mathrm{cov}(G,2)=O(\sqrt{n\log n})$ and the asymptotically sharp
estimate $O(\sqrt{n})$ relies on the Classification of Finite Simple Groups.

\bigskip

\noindent \textbf{Acknowledgement.} The author is grateful to M. Elekes for
telling him about the problem leading to this paper and for helpful
discussions.

\end{document}